\documentclass[11pt]{article}

\usepackage{epsfig}
\usepackage{graphicx}
\usepackage{color}

\newcommand{\DS}{\LARGE \renewcommand{\baselinestretch}{1.67} \normalsize}

\title{Algorithms for Deforming and Contracting Simply Connected Discrete
Closed Manifolds (I) }

\author{
Li Chen \\
Department of Computer Science and Information Technology\\
University of the District of Columbia\\
Email: lchen@udc.edu\\
}

\begin{document}
\DS

\maketitle

\abstract
{In this exploration paper, we design algorithms for
 deforming and contracting a simply connected discrete
closed manifold to a discrete sphere. Such a contraction is
a kind of shrinking or reducing process.
In our algorithms, we need to assume an ambient space for the
discrete manifold,  and this ambient space also
a simply connected discrete space in higher dimensions.

Our algorithm would work for most of cases. For some special
cases, we will make detailed analysis in the next paper. In other words,
This paper has not provided a complete proof for each case.  The algorithm
designed in this paper is in polynomial time.		
}

{\bf keyword: }{
 Discrete manifold, algorithm, contraction, deformation, sphere.}

\section{Introduction}

In 1940,   Whitehead proved that
any $C^1$ manifold $M$ has a  polyhedron $K$ and a piecewise
differentiable homeomorphism $K\to M$,
and such $K$ is a piecewise linear (PL)  manifold ~\cite{Whit40}.
Note that this polyhedron is defined as a set of triangles or simplexes.
 
In practice, no one can computationally describe a general smooth manifold or no computer can store a general smooth manifold (in any current
computing model such as Turing machines). The algorithm
must take actions on a finitely defined objects. Therefore, for the purposes of applications or realizations, we only can deal with a complex that
has finite simplexes or finite cells.

Here comes a problem, when a PL structure of a manifold is given, can we decide if it is simply connected or it is homeomorphic to a sphere? This question also has
philosophical meaning since the number of molecules in the entire universe is still finite based on our up-to-date knowledge. So it is very reasonable
that we could also ask questions on finitely defined structures.

In this exploration paper, we design the algorithm for
deforming and contracting a simply connected discrete
closed manifold to a discrete sphere (or PL-sphere).  A discrete sphere can be defined as a link of a cell in PL manifolds. It can also be defined as a PL manifold that
is homeomorphic to a sphere (in this definition, we want the deformation is very simple in discrete sense).   

More specifically, if PL manifold $M$ is simply connected,  is there any algorithm that can contract it to be a PL-sphere? We know that this problem is related to the Poincare conjecture.
It is believed that the Poincare conjecture was finally proved by Perelman in 2003. Perelman used Hamilton's Ricci flow that contains some surgical operations to some extreme cases in
his proof~\cite{KL}.

However, how do we actually deform a simply connected PL manifold to a sphere? It is also
interesting. This problem related to a constructive or algorithmic proof of this problem. 

Another related development in algorithmic aspect of this problem due to Rubinstein's work\cite{Rub}.
Rubinstein designed an algorithm to decide if a triangulated closed 3-manifold is homomorphic to the 3-sphere.
However, this algorithm is not in polynomial time. The process also needs to cut the existing tetrahedron.

In this working paper, we attempt to provide a thinking path that does not need to modify the existing simplexes or cells (no surgical operations in general).
In our algorithms, we need to assume an ambient space that is also
a simply connected discrete space in higher dimensions.		

In this paper, we will introduce a curvature-like concept called {\it the natural curviness} that will be the measure(s) for contraction or deformation. In order to prove that the algorithm always can halt,
we need to assume a statement. If this statement can be proven later, it might be true that there will be an algorithmic proof of any simply connected closed discrete manifolds
is contractible\/deformable to a $k$-sphere in discrete space. In any case, it is true that the sequence of the contractions or deformations in the process will form a homeomorphism between two discrete
manifolds: one is the original and the other is a discrete sphere that only consist of a few cells.

So the key idea of this paper is to use a deformation procedure on existing cells of an $m$-manifold $M$ to reduce the number of cells in $M$ until the number of cells to be a
small constant. See Fig. 1. As long as we can find a shorter path or sub-manifold to replace the original one, we will use the algorithm to finish the task.

In this paper, there are two basic measures for $M$: (1) Total number of $i$-cells in the $m$-manifolds $M$, and (2) Total volume covered by $M$ if we
embed $M$ into $R^n$. In our algorithm, we want either reduce the total number of of $m$-cells or total volume of the deformed $M$ in each major
step (for iteration). Because the number of $i$-cells are finite in $M$, this algorithm should be halt.
In this paper, we usually use the first one.

\begin{figure}[hbt]

\begin{center}

 \includegraphics  {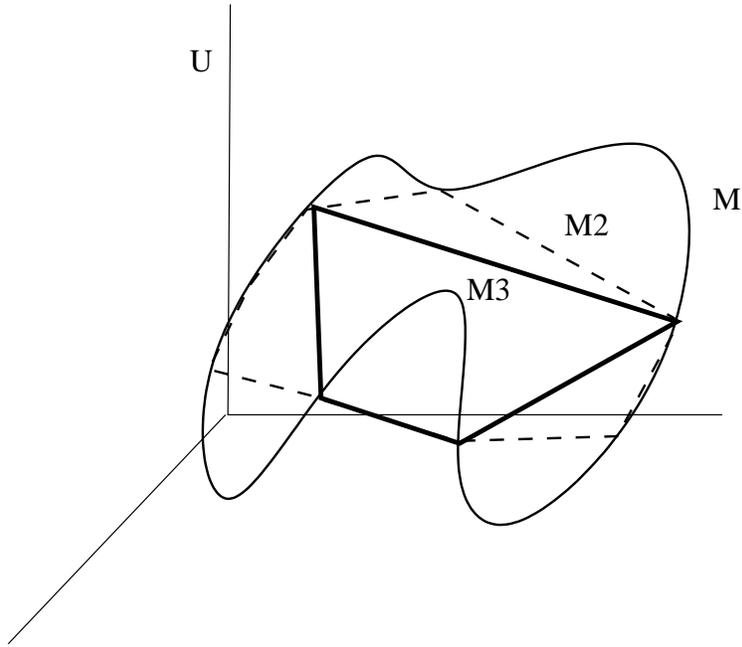} 

\end{center}
\caption{ Manifold $M$ and its deformations $M_1$ and $M_2$.  }

\end{figure}


\section{The Main Idea of the Algorithm}

In general, any smooth real $m$-dimensional manifold can be smoothly embedded in $R^{2m}$; this is called the (strong) Whitney Embedding Theorem. And any $m$-manifold with a Riemannian metric (Riemannian manifold) can be isometrically imbedded to
an $n$-Euclidean space, where $n\le c\cdot m^3$, $c$ is small constant. This is called the Nash Embedding Theorem. Therefore, we can discuss our problem in Euclidean space or a space that can be easily embedded to Euclidean space.

A simply connected (orientable) manifold $M$ in space $U$ can be viewed in Fig. 1.  If $M$ is a supper submanifold, the dimension of $M$ is smaller than the dimension of $U$ by one, in such a case, we can use Jordan's theorem to first separate the $U$ into two components. The deformation becomes the pure contraction. However, if the dimension of $U$ is much bigger than
the dimension of $M$. The contraction strategy cannot be used. So we like to use the Hamilton's idea that is to modify the local curvature of each point so that
when the curvature at each point was changed to positive or positive constant, the deformed manifold $M$ is a $m$-sphere when $m=3$. This was finally proved by Perelman in 2003.
Therefore, the Hamilton-Perelman method using Ricci flow is a local deforming method in each iteration~\cite{KL}.
 
The main idea of the Hamilton's method is to push-down the peak and to lift-up the valley geometrically. In order to make our algorithm fast, our idea is to cut out the peak and valley (
use a bridge to replace the peak or valley. )  Such a bridge is a shortest path or a minimum surface (minimum submanifold). See $M2$ and $M3$ in Fig. 1. The minimum submanifold means to count the number of $m$-cells in $M$, $M2$, $M3$ and etc.

In the process, we have the iterations too. However, such a sequence of iteration must be finite or halt in a finite time. We know that we can calculate a minimum path between two points on a curve. We also can calculate the minimum path of $i$-cells (sharing an $(i-1)$-cell) as
the distance of $i$-cells(See Appendix) between two cells.

For an arc of $M$ starting at $x$ and ending at $y$, the main idea of our algorithm is to get a type of ratio between the distance in $M$ and the distance in $U$.
(Such a distance in higher dimension is the volume.)
We will make changes at the point that has the biggest ratio or biggest reduction in terms of the number of cells in $M$.

For a curve in 2D, one can determine the sign of curvature by identifying a outward point as positive curvature the negative value to the inward point. (This is not significant in algorithm design.)

Using the same example in Fig. 1, see Fig. 2, we can see that the curve between $x$ and $y$ (or a submanifold) has the biggest ratio (in absolute sense). We want to use
the arc ${\bf a}$ to replace the original arc $A$. The key is that from original arc to the new arc, there exist a sequence of deformation. This is because $M$ and $U$ are both simply connected.

When $M$ is an $m$-submanifold, the boundary of the partial manifold (the "arc") is a closed $(m-1)$-cycle. Each $(m-1)$-cycle, $\cal C$, is homomorphic to a $(m-1)$-sphere by induction
meaning that every simply connected closed $(m-1)$- manifold ($(m-1)$-cycle) is homeomorphic to a $(m-1)$-sphere.

We know that we can get a  minimal $m$-submanifold with the boundary of $\cal C$. The meaning of the ratio is that the number of $m$-cells on $M$ bounded by $\cal C$, compare to the number of $m$-cells in the  minimal $m$-submanifold.  This minimal $m$-submanifold could have three meanings: (1) the minimum submanifold for the minimum volume in continuous sense,   (2) the minimal discrete submanifold  bounded by $\cal C$ that has the minimum number of $m$-cells,  and (3) the  discrete $m$-submanifold that passes most of  $m$-cells contained by the minimum surface in (1).

\begin{figure}[hbt]

\begin{center}

 \includegraphics  {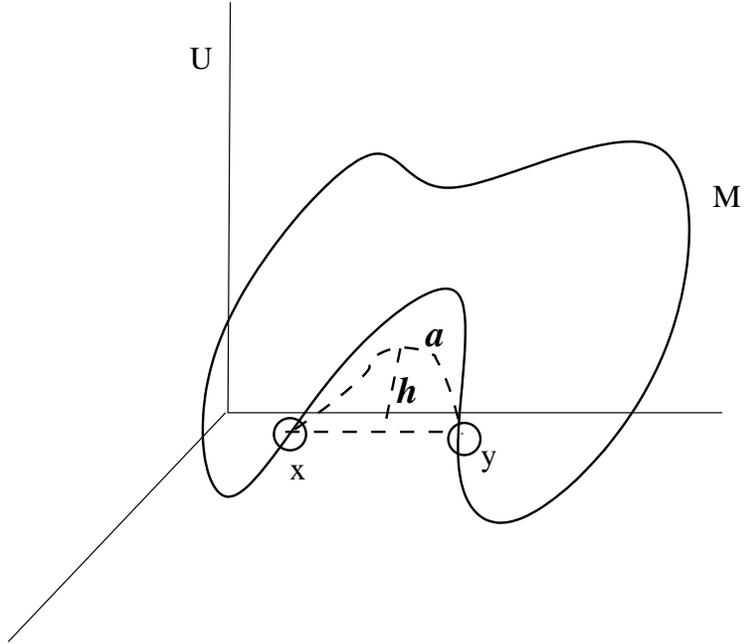} 

\end{center}
\caption{Deformation Steps}

\end{figure}

By the way, we can precalculate the distances of any two cells in $M$ and in $U$.
Ideally, for the detailed algorithm, we want to calculate the distance from each cell in ${\cal A}=M\|_C$ to the new minimum $S_C$. (The biggest one is denoted by
$h$ in Fig. 2.)    We find
the shortest paths for all, then use lofting method to get a submanifold that has less cells or less volume when closer to $S_C$, Fig. 3. (a). Here is the problem,
There are some cases in Fig. 3. (b). Some shortest paths pass ${\cal A}$ itself, in this case, there are more structures inside of the submanifold. We will split $M$ along with $C$. According to
the Jordan curve theorem or general Jordan curve theorem, this split is a complete split. So
Splitting $M$ into two components with attaching $S_C$ to both of them. To reconstruct the original manifold, we need to use connected-sum operations

\begin{figure}[hbt]

\begin{center}

 \includegraphics [width=6in]  {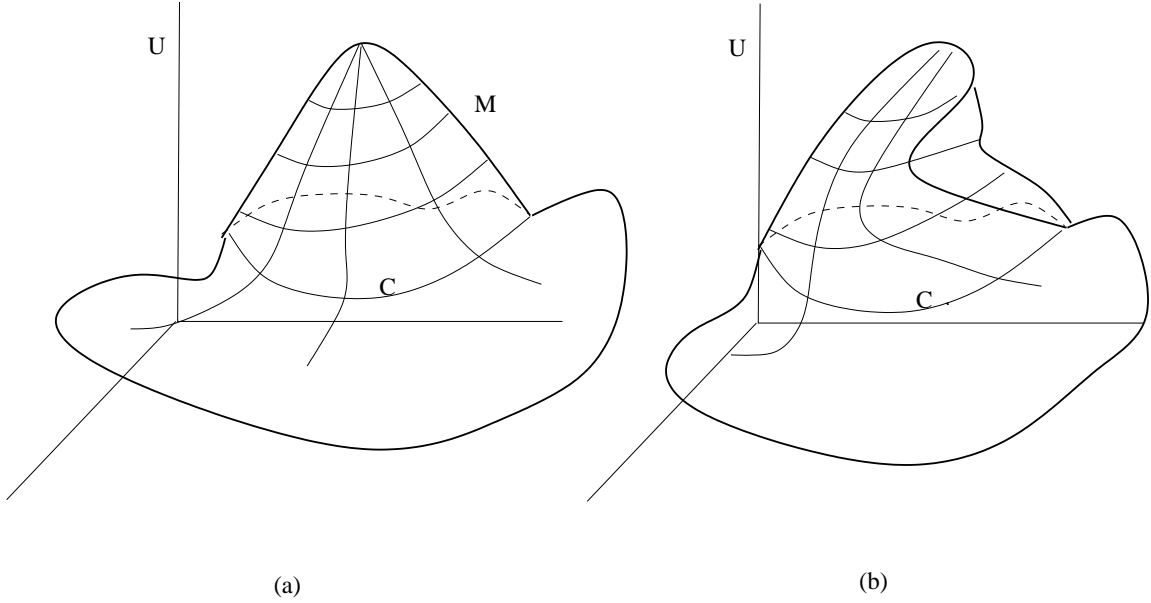} 

\end{center}
\caption{(a) A case that is deformable using our method. (b)A case should be split. }

\end{figure}

\section{Description of the Problem }

An ambient space of an object is the space that holds the object. We usually also assume
this object does not reach the boundary of the space. In this paper, we assume that the ambient space is Euclidean and it was partitioned into polyhedrons.
The modern geometry usually refers to intrinsic geometry. In 2D, this is a geometry only based on the first fundamental form of surfaces.
In intrinsic geometry, there is no ambient space to hold entire object, but it still a local Euclidean space that hold each neighborhood of a point locally.
In other words, there must be a local ambient space or a local moving frame for each local point in the object.

The computation for the connection of two local frames are not very easy to handle since we also need to sampling the data for each local points. In this paper,
we still use an ambient space to hold the object or a manifold. We have explained the rationality of this treatment by referring to the Whitney's embedding theorem and
the Nash's  embedding theorem in the first section.

The problem we will try to solve in this paper: (1) Given a discrete manifold that is simply connected and closed,  we want to deform it to be a discrete sphere. This sphere might
not perfectly rounded, but it is almost rounded. We can continue to prove that this discrete sphere is homeomorphic to a sphere.
(2) Given a discrete manifold that is simply connected and closed or a sphere, we want to reduce its size to be
only hold constant number of cells. In other words, we can easily prove it to be homomorphic to a continuous sphere.

\subsection{Smooth Manifold Triangulation in Euclidean Space}

Whitehead proved that any $C^1$ manifold $M$ has a triangulation~\cite{Whit40}. Whitney gave a method that can triangulate a differentiable manifold in Euclidean space~\cite{Whitney57}.
Whitney's method is paretically interesting since he used cubes to decompose $U$ (the ambient space) first, then  subdivide a cubical lattice called $L_0$ in $R^n$ into simplices to give a triangulation of $R^n$.
To do this is to use the center point of the each cube to make a refinement, such a refinement will generate a triangulation $L$. $L$ is called a regular subdivision of $L_0$.
In fact, one can make the lattice small enough (based on smoothness) and slightly move finite number of vertices such that the triangulation will be in general position\footnote{{\bf Definition of general position from Wikipedia}: A set of at least $d+1$ points in d-dimensional affine space (d-dimensional Euclidean space is a common example) is said to be in general linear position (or just general position) if no hyperplane contains more than d points — i.e. the points do not satisfy any more linear relations than they must. In more generality, a set containing $k$ points, for arbitrary $k$, is in general linear position if and only if no $(k-2)$-dimensional flat contains all $k$ points. }
with respect to $M$. We can generate
an $n$-dimensional simplicial subcomplex $T$ (that will be involved totally in tubular neighbourhood of $M$), The final triangulation will be a project of $T$ onto $M$.

Whitney's method provide a digitization that is particularly interesting to modern computing. In addition,  Whitney's method of triangulation is very accurate. However, it is kind of complicated.
Cairns gave a simple triangulation method for a compact closed $m$-manifold $M$ of differentiability
class $C^r$ ($r>l$) in a Euclidean space in 1961~\cite{Cairns}. His method is similar to the method of making a Voronoi diagram then a triangulation based on a dense collection of a set of points on
$M$.


\section{Algorithms of Discrete Deformation and Contraction}

We will first introduce some basic concept of algorithms and then discuss our new algorithm.
An algorithm means a sequence of finite number of instructions for completing
a task. So an  algorithm is a constructive method, but a constructive method
may not be an algorithm.

For instance, $f(x)=g(x)-x-1$ can be constructively done since we can subtract $x$ then
subtract 1 for each $x$. However, if $x$ is every point in [0,1], this constructive
procedure is not an algorithm since we could not really use a finite time process
to go through every point in [0,1].

So, strictly speaking, an algorithm must be used  for the discrete object in a finite
form.  Or at most, we can finish a task or solve a problem in countable manner. This is a
philosophical problem in mathematics. A general continuous curve is not computable in any
sense. However, it is constructible in terms of approximation. An object is
constructible does not mean that we have an algorithm to obtain that. But if an object has
an algorithm, it must be constructible.

A manifold in discrete space can be viewed as a polyhedron or a polytope, but the
edge or face does not have to be straight or flat in this paper. We call it the discretely defined cells in~\cite{Che14}. In such a way,
it is possible to generate a smooth or $C^{(k)}$ manifold using finite numbers of discrete defined cells.

Intuitively, a simply connected manifold $M$ means that every simple path in $M$ (without crossing)
from a point $p$ in $M$ and return
back to $p$ can be contracted to $p$ without crossing (each other when deforming). This type of paths is called a simple cycle.
A closed manifold means that $M$ does not have a boundary such as a circle $S^1$ and  a 2-sphere (a globe) $S^2$.

Two discrete $k$-manifolds $A$ and $B$ are called homotopy if there is a sequence of gradual variations in between.  Gradual variations
mean ``continuous'' moving from one to another one. In discrete space, ``continuous'' moving means move one unite distance at a time.
For the precise definition of discrete homotopy, please refer to~\cite{Che12}.

\subsection{Discrete Deformation Algorithms to Make Spheres}

We now assume that every manifold is orientable in this paper.
We are going to use the basic results: (1)The general Jordan curve theorem: a closed simply connected $(m-1)$-manifold $C$ on a simply closed $m$-manifold $M$ will separate
$M$ into two disconnected components. (2)The Riemann uniformization theorem: Any simply connected closed 2-manifold is homeomorphic to 2-sphere.

In this section, we assume $M$ and $C$ are both discrete and $M$ is a subset of partitioned $n$ dimensional Euclidean space $E_n$. Or $M$ is a subset of
simply connected $n$ dimensional discrete space. The homotopic mapping between two discrete $k$-submanifolds will refer to the concept in \cite{Che12}. Or
we just interpret that as the finite sequence of moves from one to another. We can make natural interpretation. We also assume that $M$ contains
finite cells and each cell has relative the same volume in each dimension (with respect to a constant).

Again, as we discussed in the previous section,
the main idea of the algorithm in this paper for making spheres is to push the ``peak'' and to lift the ``valley.''  This idea was from Hamilton's plan for proving
Poincare conjecture. So finding the peak or valley with largest ``slope''  will be the key to our algorithm. This also can be implemented by using
a flatter one with much less cells to replace the existing ``peak'' or ``valley.''

Here comes a problem, for a curve in 2D, a ``peak'' is a outward point and a ``valley'' is a inward point. However, in 2D surface,  the total (Gaussian) curvature
has different meaning to the mean curvature,  the ``valley'' points also have positive Gaussian curvature. The negative Gaussian curvature points are hyperbolic points
that are most flat (closer to mean curvature to be zero, sometimes). So using (Gaussian) curvatures to be scalar curvature for measuring the curviness may not be very smart.

Therefore, we can just use the natural curviness that is only care of the number of cells in the manifold. We want to reduce them to get a small constant number.
We will not treat hyperbolic points in the algorithm and let the process of modifying ``peaks'' and ``valleys'' automatically change the shape of  hyperbolic points in
the finite discrete object. This is based on the Gauss-Bonnet theorem for genus zero manifolds where we always can find positive outward points.

Let us introduce some formal definitions for this purpose. In $M$, an arc $A$ is an $m$-dimensional submanifold of $M$ with a boundary usually called $C$, $C$ is an discrete
$(m-1)$-manifold. A chord $B$ regarding to $C$ is a minimum volume (length or area etc)
$m$-submanifold (line or surface, as a ``base'' in $U$) with the same boundary $C$.  The hight $h$ is the shortest path length for any point in $A$ to $B$ in space $U$. In discrete
cases, from $A$ to $B$ we only need $h$ times deformation or we only need to insert $h-1$ gradually varied surfaces\/submanifold in between (if there is no holes in between). Such a deformation
generate a type of homomorphism in discrete sense in ambient space $U$.


We now define a function for measuring the ``slope.''. First, we start with curves, we assume $M$ is a curve in $U=R^n$ the distance from a point $p$ to $q$  in $M$ is denoted
as $D(p,q)=d_{M}(p,q)$ and  $d(p,q) = d_{U}(p,q)$ indicates the shortest path in $U$ or in the ambient space. $d_{M}(p,q)$ is just   the number of points from $p$ to $q$ in $M$ (or subtract by one).
There are total of four measures that are indicating the natural curviness of the arc between $p$ and $q$ in M.

\[r(p,q)=\frac {d_{M}(p,q)}{d(p,q)}\]

\[r^{(1)}(p,q)= d_{M}(p,q)-d(p,q)\]

\[r^{(2)}(p,q)= h \]

\[r^{(3)}(p,q)= \frac {h}{d(p,q)} \]

\noindent When we deal with higher dimensional manifolds, $d_{M}(.)$ is the volume of the arc $A$, and $d(.)$ is the volume of the minimum manifold. $r^{(3)}(p,q)= \frac {h}{d(p,q)}$ will
be the $r^{(3)}(C)= \frac {h}{d(B)}$, $d(B)$ is the diameter or average diameter of $B$. All purposes are the reduction of the number of cells in $M$. Here volume means primarily the number of cells.
However for some extreme cases, number of cells may be the same . So we need to use volumes in addition to that.  Now we define , $N(A)$ as the number of cells in highest dimension in $A$.
$vol(A)$ is the volume when embedding to a manifold with a metric.  In general, we use $Vol(A)$ to represent $N(A)$ with consideration of $vol(A)$.

Formally, $Vol(A)>Vol(B)$ means that (1) $N(A)>N(B)$, or (2) $N(A)=N(B)$ but  $vol(A)>vol(B)$.
 
If $h$ is 1, we know that $A$ is gradually varied to $B$. There is a possibility that $N(A)=N(B)$ but $vol(A)>vol(B)$. So we can still use $B$ to replace $A$.  This is
only valid to be used when $N(A)=N(B)$.

Now we can define the natural curviness at a point $x$ in $M$? We can use the radius $\gamma$ to help us to measure.  Since now $M$ is a curve (1-manifold, or other), we can find $p$ and $q$ such that
$D(p,x)=\gamma$ and $D(q,x)=\gamma$.

\[r_{\gamma}(x) = r(p,q). \]

\noindent where $r(p,q)$ can be any of $r^{(i)}(p,q)$ defined above with fixed $i$.

\noindent Since $M$ is close, when $\gamma$ is big,  $r_{\gamma}(x)$ could be very big since the path corresponding to $d_{M}(p,q)$ may cover all points in $M$. So, we may require such $arc (p,q)$ containing less than half of points in $M$. So we will start with small radius $\gamma$ to check which peak or valley  will be picked first for processing.

If $M$ is a surface or high dimensional $k$-manifold, we can still use the similar idea. We can use a closed $(k-1)$-manifold $C$. The distance from each point ($(k-1)$-cell) in $C$ to $x$ will be $\gamma$ (or very close to $\gamma$, does not have to be exact as the value $\gamma$ in computation).
In discrete space, $\gamma$ is an integer.  We know that $C$ will split $M$ into two components. We like to see that one contains less than half $k$-cells of $M$.
we can use $m(C)$ as the minimum surface ($(k-1)$-manifold) bounded by $C$. $m(C)$ is a discrete $(k-1)$-manifold that contains minimum number of $k$-cells (or $Vol(.)$) in the ambient space. Some or many of
the cells might not in $M$.

Formally, when given a radius $\gamma$ , we need to find a boundary of $k$-cell $p$ in $M$ such that

\[ B(x,\gamma)=\{ p | D(p,x) \le \gamma \} . \]

\noindent Ideally, $B(x,\gamma)$ is a ball. But in discrete case, this may not be true. It may be dependent on a decomposition. However, we can get a closed $C$ as the boundary that is a discrete $(k-1)$-manifold,
 the smallest, and containing $B(x,\gamma)$.   We can also get a closed $C_0$ as the boundary that is a discrete $(k-1)$-manifold,
 the biggest, and is contained by $B(x,\gamma)$.  If the decomposition is a triangulation, we believe that $C_0=C$, or just off by 1 in radius.

Algorithmically, we need to proof that this set forms a subset of manifold $M$ with a cycle boundary. If not, we need to find a boundary cycle $C$ that is a good fit to this set. We mean that maximizing the number of
cells regarding to  $B(x,\gamma)$. This closed curve\/boundary $C$ may not be unique, but the number of cells in $M$ included the submanifold with boundary  $C$ is unique. The set is denoted as ${B_{\gamma}(C)}$.  We use $C(x)$ to denote this cycle.  Then we find
$m$-submanifold with the minimum cells bounded by $C(x)$.  We denoted this submanifold as $m(C(x))$, a minimum surface.

For example we can use

\[r_{\gamma}(x) =  \frac {Vol(B_{\gamma}(C(x))}{Vol(m(C(x)))}  \]

\noindent or,

\[r_{\gamma}(x) =  Vol(B_{\gamma}(C(x))-Vol(m(C(x)))  \]

\noindent as the measure for the natural curviness. Note that $B_{\gamma}(C(x))$ and $m(C(x))$ form a closed manifold, ideally. We will discuss the pathological cases later.
In addition,  if it costs too much time to get such a best $C(x)$,  we can just use an approximation as long as we can reduce the volume. This should be good enough in practice.


\subsection{Finding Sign of an Arc in $M$}

In this paper, what we need to do is to use algorithms to reduce the size of the manifolds to be the one that contains constant numbers of $m$-cells.
The method suggested by Hamilton is to change the curvature of each point to be positive curvature. After that, it can be done to make every point to have
a constant curvature. When a manifold has a positive constant curvature, then we can see it is a sphere.

This idea is based on the fact that the integration of Gaussian curvature of each point is a positive number for a simply connected closed manifold. So it is obvious that
objective can be made in such a way.

In discrete cases and computations, it is hard to find Gaussian curvature  without errors. Even though, it is possible to do this theoretically. But for practice,
it may cause some other problem such as precisions.
 
To determine the sign of curvature at $x$ may cause some problems for some manifolds.
 
For a simple example, see Fig. 4. If we define clockwise as our "moving"  direction then the original part in the (future) separated 1-cycle, is clockwise
then we will say that the sign is positive. We find the peak and we need to reduce the height of the peak. Otherwise, we need to lift the valley. Arc $a$ is
on a plane that can be outward arc or inward arc depending on how we make a filling bounded face (an $m+1$-manifold with boundary $M$).

\begin{figure}[hbt]

\begin{center}

 \includegraphics  {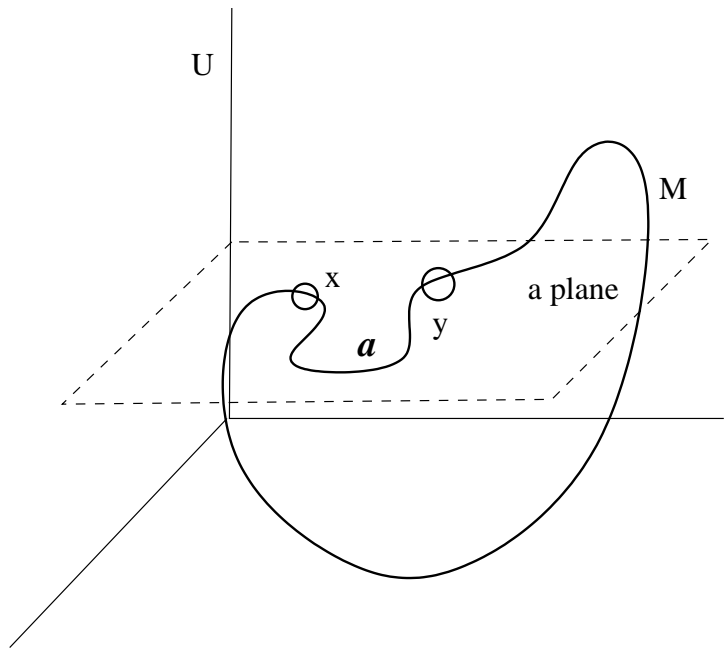} 

\end{center}
\caption{ The sign of curvatures: (a) Negative , (b) positive }

\end{figure}

For the Gaussian curvatures of surfaces, the inward and outward points do not determine the sign of curvature. However, we will lift inward points and push the outward points.
As long as we can reduce the volume, we will do it.
If an inward part contains an outward part. This is complex component. We just want to cut it and
treat it and than use connected sum to make it back as needed for homomorphic mapping.

Any arc in $(m+2)$-space, can be either positive or negative curvature (or zero). It depends on how we make $M$ to be the boundary of $(m+1)$-manifold. This
is because we can make inward or outward arc by making different $m+1$-cells.   We can make it as always positive curvature in space. So we can just use the Hamilton's theorem on positive curvatures.
This example can be an argument that the intrinsic geometry is enough for
describing a geometric object. See Fig. 5. even though $M$ is on the a plane, in Fig. 5 (a), Arc $a$ is inward arc and the sign is negative. But in Fig. (b),  Arc $a$ is outward arc and the sign is positive.

Only for the case $M$ in $m+1$-space, we will have the  definitive negative or positive curvatures for each arc. The filling is unique due to the Jordan Curve Theorem.

More discussions are below.
Even though, there is an inward (negative curvature) arc. We can cut the inward part by filling such a part into the space to make it none-negative.

In discrete space, this can be done using the algorithm described in this paper.
Therefore, local curviness cannot determine the sign of the curvature. The sign (positive or negative ) will be determined by an filling (an interpretation ) on the closed manifold $M$ that is
a $(m+1)$-manifold. In other words, except zero, the sign of the curvature of an arc (or at a point) may change based on a filling.

\begin{figure}[hbt]

\begin{center}

 \includegraphics [width=4.5in] {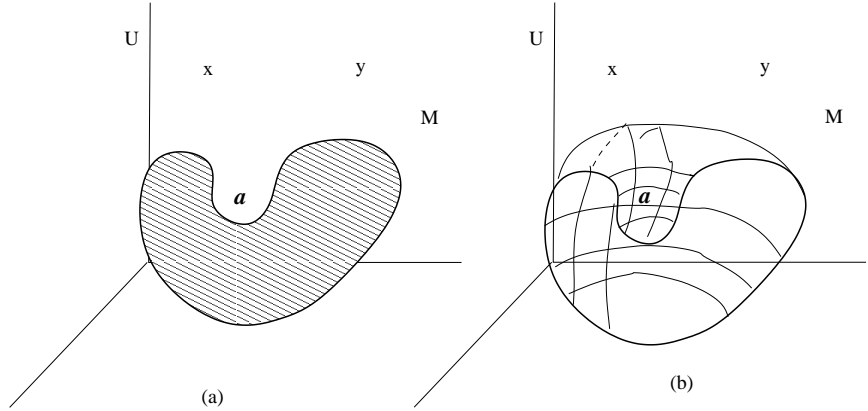} 

\end{center}
\caption{ One boundary with two interpretations: (a) Negative curvature, and (b)positive curvature }

\end{figure}

\subsection{Discussion About the Discrete Distances in Manifolds}

We said that we will calculate the distance for every pair of points in the manifold. How to select a radius $\gamma$? We can first select the largest distance pair $p$ and $q$, $D(p,q)=d$. Let $\gamma=[d/2]$, Then $B(p,\gamma)$ and $B(q,\gamma)$ will have some common cells in their boundaries. If we can find a $C$ pass most of cells in their boundaries, we may get a separated two manifolds.
We can see : (a) If $m(C)$ does not contain any cells in $M$ other than $C$, so we think that $C$ split $M$ as two almost equal components. When $Vol(m(C))$ is bigger than the volume of any of the components, we can see that
we want to find examine two other pairs with the same distance. (2) When $m(C)$ contains some cells in $M$ other than $C$, we can let $\gamma:=\gamma/2$ or  $\gamma:=\gamma-1$ to search for a radius. So
we can always find one $B(p,\gamma)$ that is a peak or valley.

How do we find an outward points in $M$? we can find the $d(x,y)$ in $U$ for all $x$,$y$ in $M$. There must be a pair that have the largest distance $p$ and $q$ in $M$, denoted as $D$, the (discrete) diameter. So the neighborhood of  $p$ or $q$ must be a outward (peak) point area.  If for every $x$ in $M$, there is a $y$ in $M$ such that $d(x,y)=D$ (or almost equal to $D$), then $M$ is a sphere. This is also an advantage of the discrete method. We
do not have to use curvature to exam the roundness. We can use algorithm to check one by one.


\subsection{The Main Algorithm}

The key part of the algorithm is to modify a section or part of $M$ if this part $X$ (to be centered at point $x$) is not a minimum submanifold $m(X)$ (in discrete sense). Such a modification is
to reduce the volume of $X$ to approximate the minimum step-by-step. It is
done by its gradual variation (the simplest deformation) on its neighboring envelop $S(X)$.  In order to guarantee such a $S(X)$ is a submanifold in $U$. $X$ need to be flat(local fateness).

Using minimum surface or submanifold to replace a positive curvature parts (in 3D ) will add more positive peak and valleys. This is because a minimum surface is usually with the negative total curvature.
So the philosophy of this process of removing peaks and valleys is very good: while reduce the total volume of the manifold M, we are adding more positive curvature points. In other words,
the number of peaks are increased.
So on so fourth, we can get a smallest (irreducible) sphere-like discrete manifold that only contains a few $m$-cells.

While doing the big curviness point removing,
we need assume such a $m(X)$ does not contain any cell in $M$ except at the boundary of $m(X)$. In such a case, we need to select a part $Y$ that is in $X$ to make the deformation first. Such $Y$ is always
exist due to the fact of finiteness of $M$.

In order to find such a deformation fast, we need another concept called lofted surfaces (submanifolds). Find a circle (cycle) on $M$ centered at $x$, $C_i(x)$, every point on this cycle has the same graph-distance to
$x$ that is $i$. Make a minimum surface $m_i$ with the boundary $C_i(x)$. Usually $m_i$ dose not contain any cell in $M$ (any type of cells) except on $C_i(x)$.  Now, if we have such a sequence $m_{i}$, the
set bounded by $X$ and $m(X)$ is called a semi-convex. If not, $X$ will be still replaced by $m(X)$, but we will do special treatment for the closed manifold $X\cup m(X)$. Actually, we will do recursive process on
this new manifold. As we talked that $X$ usually contain less than 1/2 of $m$-cells in $M$. But any way, $m(X)$ is minimum so $X\cup m(X)$ is smaller then $M$ in terms of size.  See Fig. 6.

\begin{figure}[hbt]

\begin{center}

 \includegraphics [width=6in] {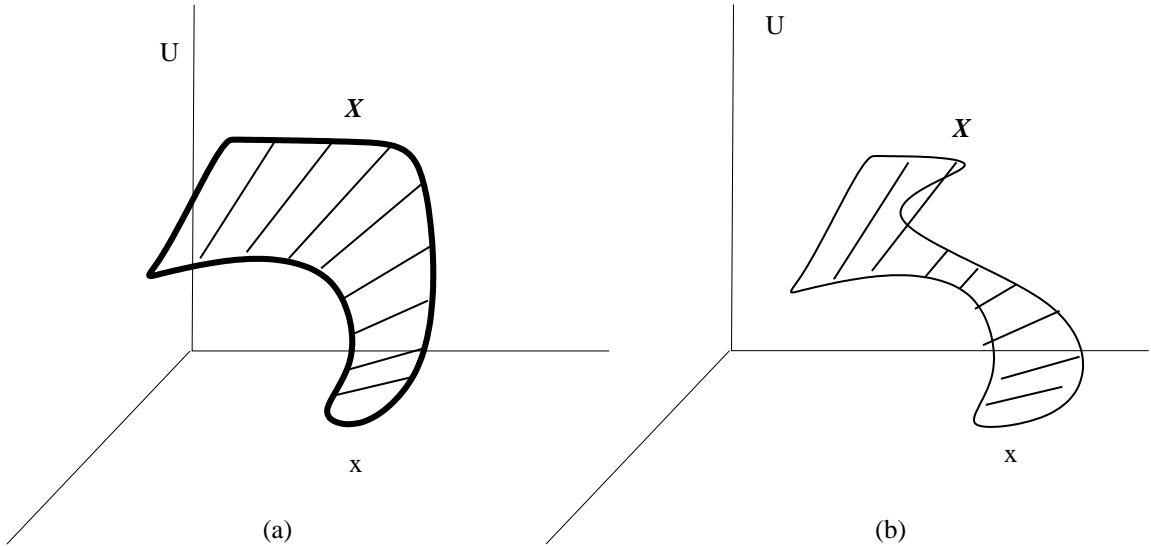} 

\end{center}
\caption{Sequences of minimum circles: (a) the regular case, and (b)the irregular case   }

\end{figure}

As we can see that, if we have such a sequence of $m_i$, we can find a modification (deformation) a long with $m_i$.  As long as we can find a $m_i$ that is not original part of $X$, the size is smaller (not equal).
So we can use the $m_i$ to replace the old one.  The curviness is improved at least by one cell toward to the minimum. Even though the volume of $m_i$ is not decreased monotonically, as long as each pair of $m_i$ and $m_{i+1}$
are gradually varied, so the process of using $m(X)$ to replace $X$ in $M$ also works.

We can see that $M$ is always reduce in it size in removing a peak or a valley. If the size cannot be reduced, then $M$ is already cannot be cut. Ideally, such $M$ contains all $m$-cells where each cell will contain the certain point $q$ . This unique point
$q$ is the center of the $(m+1)$-sphere in discrete sense. But such a $q$ is not in $M$.  (We can also use one more cut to make $M$ is a $(m+1)$-cell in discrete sense. then contract this $(m+1)$-cell to be one point.)

If we do not assume $M$ is simply connected, then at some point, the set bounded by $X$ and $m(X)$  is not semi-convex.   So this algorithm can determine if M is simply connected.  For instance, for a tours, we will
always intersect points other than boundary of $m_i$ in $M$.

{\bf Algorithm 3.1, The main algorithm:}  In this algorithm,  we will push the highest positive curviness part by deformation, if we can not deform this part meaning it contains complex structure, we will cut it out to do recursive analysis.

\begin{description}
	
\item [Step 1] Make $M$ to be a local flat $m$-manifold. This is done by locally small modification of $M$ in $U$. Later on, each removing or reduction will not generate a non-flat area. Because we will always use $m(X)$, the minimum submanifold to replace the original $X$. $m(X)$ is locally flat.

\item [Step 2] Calculate the discrete natural curviness $r$ for each radius $\gamma$ at every point $x$ on $M$. (Find the minimum discrete surface $m(C(x))$ for each $B(x,\gamma)$.  $C(x)$ in  $M$ is a discrete $(m-1)$-closed manifold and it is the best approximation of all points whose distance to $x$ is the radius $\gamma$ on $M$. Choose the set so that $m(C(x))$ does not contain any cell in $M$ except in $C(x)$. Then choose the set so that
     $Vol(m(C(x)))$ is smaller than each of the two components separated by $C(x)$ of $M$. This set is our valid set called  $\Gamma$. Any of the member can be used in our later steps.  If  $\Gamma=\emptyset$,
     there must be a point $o$ in $U$, such that every $m$-cell in $M$ will interest with an $m$-cell in $U$ containing $o$. This is a discrete sphere. It is homotopic to the point $o$.)

\item[Step 3] Find a point $x$ the maximum $r$ starting at $\gamma=1/4$ radius (radius=$D/2$ of $M$) (then reduce the radius, from small radius to bigger is fine too). Let $X$ is such a disk with regarding to $M$. (Do the modification or deformation above).

\item[Step 4]  In discrete case, the boundary of $X$ might not be $(m-1)$-cycle (discussed previously), we find $C$ that is a $(m-1)$-cycle and best fit to the boundary of $X$. By the Jordan curve theorem, $C$ split $M$ into
two components, one of which contains $X$ with minimum number of other cells not in $X$. There are two kinds of minimum surfaces\/manifolds regarding to $C$ now, one is the continuous minimum surface $min(C)$ and the discrete
minimum digital surface $m(C)$. (There are may be multiple $m(C)$'s) We select $m(C)$ that contains most of cells in (intersecting with) $min(C)$.

\item[Step 5] Make a discrete minimum surface based on $min(C)$, we still use $m(C)$ to represent that is a discrete surfaces containing most of cells in $min(C)$ with the minimum $m$-cells. Containing most of cells in $min(C)$ is the primary needs. (For every $x$ and $\gamma$, if $m(C)$ only contains one cell that is not in $M$, that means $M$ is already a sphere and such a cell is the center of the sphere.)

\item[Step 6] If  $m(C)$ contains other cell in $X$ other than $C$, do not continue. It means that $X$ contains an inner part that interests with the minimum surface. We did not find the right $X$. Change to another $X$
with the same radius go to Step 3. If no more such $X$ in the same radius; reduce the radius to find new $X$ in Step 3.

\item[Step 7] Continue Step 5, we find an $m$-cycle $X\cup m(C)$. We now want to deform $X$ to be $m(C)$. There are two ways to do it. We present a way here that is called lofted-circles. Let $C_i=\{y | D(x,y)=i\}$ and $i=1,\cdots,\gamma$. $C_i$ is called the lofted circle. We still use $C_i$ to represent a closest $(m-1)$-cycle to $C_i$. $X$ restricted by $C_i$ is denoted by $X_i$.  Let $m_{i}$ be the minimum surface for $C_i$.
    If a cell in $X$ interests with $m_i-C_i$ , we stop to continue here for deforming $X$ to $m(C)$. We cut $X$ out and  use $m(C)$ to replace $X$. We will use this algorithm recessively to the manifold $X\cup m(C)$ later.
    It means in this case, we split $M$ into two manifolds, we can later use connected sum to glue them together to get the homeomorphism.

\item[Step 8] If  $X\cup m(C)$ is a semi-convex based on $m_i$ calculated in Step 7. (No cell in $X$ interests with $m_i-C_i$.) Now we can start the deformation from $X$ to $m(C)$.

\item[Step 9] {\it Sub-Algorithm A}: We need to make sure or prove now, such a lofted surface will bring a sequence of discrete surface for certain. These surfaces are mostly gradual varied. In some occasions, we need to insert some deformation surfaces in between $m_{i-1}$ to $m_{i}$ .     What we can do here is to use a special local modification method to find the subsequence inserted in between  $m_{i-1}$ to $m_{i}$. (Note that We use the same idea of meeting the most of cells in $m_i$ and the search is at the near $S(X_i)$, the star or neighborhood of $X_i$.) We will give more detailed steps of this Sub-algorithm later.

\item[Step 10] We get a homomorphism mapping by deforming the higher discrete curviness part to a minimum surface. This process will add more positive curviness parts in theory.  This process will at least reduce a $m$-cell from $M$ unless the discrete curviness is the same for all points. Or a discrete sphere is arrived, so there is no minimum surface can pass a $m$-cell that is not in $M$ ($M$ is a $(m+1)$ cell).  Or every minimum surface
    will contain the same point not in $M$ that is the center of a sphere. In this case, $M$ only contains very limited cells. and each $m$-cell to this point has the constant distant 1.

\item[Step 11] Repeat this process from Step 2.   We will have a discrete sphere.

\item[Step 12] When we deal with discrete minimum surfaces , the surface must at least contain a cell that is not on $M$.

\end{description}

In fact, as long as the discrete natural curviness exist for any point $x$ for any $\gamma$, we can continue our process for reducing the number of cells in $M$. This reduction process will be halt until there is no minimum surface that contains a cell what is
not in $M$. Or the minimum surface only contains one cell that is not in $M$. If for every part, we have the same cell (or point) $O$ in all minimum surfaces.  Such an $M$ only contains constant number of $m$-cells that is centered
at this special cell $O$.

We will give the detailed steps of {\it Sub-Algorithm A} in the next paper in the near future.

\section{Summery}

The purpose of this paper is to explore the constructive and algorithmic method for deformation and contraction. We try to not use surgical operations in topology. 
It is our hope that in our approach, we can avoid some singularities. However, the difficulty is that it is hard to describe a discrete sphere and find the perfect 
match between discrete minimum and continuous minimum. There are still much work to do in the future.

\section{Appendix: Some Concepts in Manifolds and Discrete Manifolds }

The concepts of this paper are in \cite{Che14}. We also use some concepts from the following two papers:
L. Chen    A Concise Proof of Discrete Jordan Curve Theorem, {\tt http://arxiv.org/abs/1411.4621} and
L. Chen and S. Krantz,  A Discrete Proof of The General Jordan-Schoenflies Theorem,
\newline {\tt http://arxiv.org/abs/1504.05263}.

{\it
A discrete space is a graph $G$ having an associated structure. We always assume that $G$ is finite, meaning that $G$ contains only a
finite number of vertices.  Specifically, ${\cal C}_2$ is the set of all minimal cycles representing all possible 2-cells;
$U_2$ is a subset of  ${\cal C}_2$.   Inductively, ${\cal C}_3$  is the set of all minimal 2-cycles made by $U_2$.   $U_3$ is a subset of  ${\cal C}_3$ .
Therefore $\langle G,U_2,U_3,\cdots,U_k \rangle$ is a discrete space. We can see that a simplicial complex is a discrete space.  For computational purposes,  we want to require
that each element in $U_i$ can be embedded into a Hausdorff space or Euclidean space using a polynomial time algorithm （or an efficient constructive method). And such a mapping will
be a homeomorphism to an $i$-disk with the internal area of the cell corresponding to an $i$-ball that can be determined also in polynomial time.
Another thing we need to point out here
is that $u\cap v$ in   $\langle G,U_2,U_3,\cdots,U_k \rangle$ must be connected. In most cases,  $u\cap v$ is a single $i$-cell in $U_i$ or empty.
In general,  $u\cap v$ is homeomorphic to an  $i$-cell or empty.  In \cite{Che04,Che14}, we used connected and regular points to define this idea for algorithmic purposes. This is because
the concept of homeomorphism is difficult for calculation. Now we request:  that $u\cap v$ is
homeomorphic to an  $i$-cell in polynomially computable time. We also would like to restrict that idea to decide
if an $i$-cycle is a minimal cycle or an $(i+1)$-cell is also   polynomial time computable. As an example,  a polyhedron partition usually can be done in polynomial time in computational geometry.

In our definition of discrete space (a special case of one such is PL space, meaning that our definition is more strict),
a $k$-cell is a minimal closed $(k-1)$-cycle. A minimal closed $(k-1)$-cycle might not be a $k$-cell in general discrete space since it is
dependent on whether the inner part of the cell is defined in the complex or not.  We view that a $1$-cycle is a closed simple path
that is homeomorphic to a $1$-sphere. So  a $(k-1)$-cycle
is  homeomorphic to a $(k-1)$-sphere. The boundary of a $k$-cell is a $(k-1)$-cycle.


We also need another concept about regular manifolds. A regular $k$-manifold $M$ must have the following properties:
(1) Any two $k$-cells must be $(k-1)$-connected, (2) any $(k-1)$-cell must be contained in one or two
$k$-cells, (3) $M$ does not contain any $(k+1)$-cells, and (4) for any point $p$ in $M$, the neighborhood of $p$ in $M$,
denoted by $S(p)$, must be $(k-1)$-connected
in $S(M)$.

In the theory of intersection homology or PL topology~ \cite{GM},  (or as we have proved in \cite{Che13}), the neighborhood of $x$
(containing all cells that contains $x$) $S(x)$ is called the {\it star} of $x$. Note that $S(x) \setminus \{x\}$ is called the {\it link}.
Now we have: If $K$ is a piecewise linear $k$-manifold, then the link $S(x) \setminus \{x\}$ is a piecewise
linear $(k-1)$-sphere.
So we will also write ${\rm Star}(x)$ as $S(x)$ and ${\rm Link}(x)={\rm Star}(x)-\{x\}$.  In general, we can
define ${\rm Star}({\rm arc}) = \cup_{x\in {\rm arc}} {{\rm Star}(x)}$. So ${\rm Link}(\rm arc) = {\rm Star}({\rm arc}) - \{arc\}$.
${\rm Star}({\rm arc})$ is the envelope (or a type of closure) of ${\rm arc}$.

We also know that, if any $(k-1)$-cell is contained by two $k$-cells in a $k$-manifold $M$, then $M$ is closed.

In a graph, we refer to the distance as the length of the shortest path between two vertices.
The concept of {\it graph-distance} in this paper is the edge distance, meaning how many edges are needed from one vertex to another. We usually use the length of the
shortest path in between two vertices to represent the distance in graphs. In order to distinguish from the distance in Euclidean space,
we use graph-distance to represent lengths in graphs in this paper.

Therefore graph-distance is edge-distance or 1-cell-distance. It means how many 1-cells are needed to travel from $x$ to $y$.  We can generalize
this idea to define 2-cell-distance by counting how many 2-cells are needed from a point (vertex) $x$ to point $y$.  In other words,
2-cell-distance is the length of the shortest path of 2-cells that contains $x$ and $y$. In this path, each adjacent pair of
2-cells shares a 1-cell. (This path is 1-connected.)

We can define $d^{(k)}(x,y)$, the $k$-cell-distance from $x$ to $y$, as the length of the shortest path of (or the minimum of number of $k$-cells in such a sequence)
where each adjacent pair of two $k$-cells shares a $(k-1)$-cell. (This path is $(k-1)$-connected.)

We can see that $d^{(1)}(x,y)$ is the edge-distance or graph-distance. We write $d(x,y) = d^{(1)}(x,y)$

(We can also define $d^{(k)}_i(x,y)$) to be a $k$-cell path that is $i$-connected. However, we do not need to use such a concept in this paper. )}

\end{document}